\documentclass[11pt]{amsart} 
\title[Operator Space Projective Tensor Product]{Operator Space Projective Tensor Product: Embedding into Second dual and ideal structure}

\usepackage[all]{xy}
\usepackage{amsmath,amsfonts,amssymb,amsthm,a4,hyperref}
\newtheorem{thm}{\sc Theorem}[section]
\newtheorem{cor}[thm]{\sc Corollary}
\newtheorem{prop}[thm]{\sc Proposition}

\newtheorem{rem}[thm]{\sc Remark}
\newtheorem{lem}[thm]{\sc Lemma}
\newenvironment{pf}{\noindent {\sc Proof:}}{\hfill $\Box$}

\newcommand{\oop}{\widehat\otimes}
\newcommand{\oh}{\otimes_h}
\newcommand{\omin}{\otimes_{\min}}
\newcommand{\obp}{\otimes_\gamma}

\newcommand{\ra}{\rightarrow}

\newcommand{\C}{\mathbb{C}}

\newcommand{\B}{\mathcal{B}}
\newcommand{\K}{\mathcal{K}}  
\newcommand{\CB}{\mathcal{CB}}
\newcommand{\JCB}{\mathcal{JCB}}
\newcommand{\N}{\mathbb{N}}
\newcommand{\M}{\mathbb{M}}
\begin{document}

\author[R. Jain]{Ranjana Jain}
\address{Department of Mathematics\\ Lady Shri Ram College for Women\\ New Delhi\\ \hspace*{3mm} India.}
%\curraddr{Department of Mathematics\\ K.U. Leuven\\ Belgium.} 
\email{ranjanaj\_81@rediffmail.com}
\thanks{The research of the first author was supported in part by the Research Foundation Flanders (FWO), under the Research Programme G.0231.07}
\author[A. Kumar]{Ajay Kumar}
\address{Department of Mathematics\\ University of Delhi\\ Delhi\\ India.}
\email{akumar@maths.du.ac.in}
%---------------------------------------------------------------------------------


\begin{abstract} 
   We prove that for operator spaces $V$ and $W$, the operator space $V^{**}\oh W^{**}$ can be completely isometrically embedded into $(V\oh W)^{**}$, $\oh$ being the Haagerup tensor product. It is also shown that, for exact operator spaces $V$ and $W$, a jointly completely bounded bilinear form on $V\times W$ can be extended uniquely to a separately $w^*$-continuous jointly completely bounded bilinear form on $ V^{**}\times W^{**}$. This paves the way to obtain a canonical embedding of $V^{**}\oop W^{**}$ into $(V\oop W)^{**}$ with a continuous inverse, where $\oop$ is the operator space projective tensor product. Further, for $C^*$-algebras $A$ and $B$, we study the (closed) ideal structure of $A\widehat{\otimes}B$, which, in particular, determines the lattice of closed ideals of $\B(H)\widehat{\otimes} \B(H)$ completely.
\end{abstract}

\maketitle

\section{ Introduction }
The operator space projective tensor product serves as an analogue to the Banach space projective tensor product in the category of operator spaces. It is used to linearize the jointly (matricially) completely bounded bilinear maps in the same way as the Banach space projective tensor product linearizes the bounded bilinear maps. If $E$ and $F$ are operator spaces, then their {\it operator space projective tensor product}, denoted by $E\widehat{\otimes}F$, is the completion of the algebraic tensor product $E\otimes F$ under the norm
$$
\|u\|_{_{\wedge}} = \inf
\{\|\alpha\|\|v\|\|w\| \|\beta\| : u = \alpha (v \otimes
w) \beta \}, \, u \in M_n(E\otimes F),
$$
where the infimum runs over arbitrary decompositions with $ v\in M_p(E),w\in M_q(F), \alpha\in \M_{n,pq}, \beta\in \M_{pq,n}$ and $ p, q\in \N $ arbitrary, $\M_{k,l}$ being the space of $k\times l$ matrices over $\mathbb C$. The theory of operator space tensor products was developed independently by Blecher and Paulsen \cite{bp}, and Effros and Ruan \cite{er3,er4}.

For $C^*$-algebras $A$ and $B$, it is known that $A^{**} \oh B^{**}$ can be completely isometrically embedded into the bidual $(A \oh B)^{**}$ \cite[Theorem 4.1]{kumsi}. In Section 2, we prove that the same is true, in general, for operator spaces. Haagerup, in \cite{haag}, proved that a bounded bilinear form on $A\times B$  has a unique norm preserving extension to a separately normal bounded bilinear form on $ A^{**}\times B^{**}$. We prove an analogous result for jointly completely bounded bilinear form on exact operator spaces and for arbitrary $C^*$-algebras. Using the above extension result of Haagerup, Kumar and Sinclair \cite{kumsi} proved that for $C^*$-algebras $A$ and $B$, there is a canonical bi-continuous embedding of $A^{**} \otimes_{\gamma} B^{**}$ into $(A\otimes_{\gamma}B)^{**}$, where $\otimes_{\gamma}$ denotes the Banach space projective tensor product. Its analogue for the operator space projective tensor product was left open. In this paper, we present an affirmative answer in the cases of exact operator spaces and arbitrary $C^*$-algebras. As an application of this result, we re-establish the fact that for $C^*$-algebras $A$ and $B$, the Haagerup norm and the operator space projective norms are equivalent on $A\otimes B$ if and only if $A$ and $B$ are subhomogenous, which was already proved in \cite{kumsi}.  It must be remarked here that Sanchez and Garcia \cite{san} also studied the relationship between the tensor product of the biduals and the bidual of the tensor product for Banach space projective tensor norm. In particular, they proved that for a Banach space $X$ of type 2 such that $X^*$ is of cotype 2, the embedding $X^{**}\obp X^{**}\hookrightarrow  (X\obp X)^{**}$ is bi-continuous.

The closed ideal structure of $C^\ast$-algebras $A$ and $B$ for $A\otimes_{h}B, \, A\otimes_{\text{min}}B$ and $A\otimes_{\text{max}}B$ has been investigated by Allen, Sinclair and Smith \cite{ass}, Archbold, Kaniuth, Schlichting and Somerset \cite{arch}, Takesaki \cite{take} and Wassermann \cite{wass} respectively. For the commutative case, the closed ideals of $A\otimes_{\gamma}B$ have been discussed in Graham and McGehee \cite{gmg}. However, the analysis of the (closed) ideal structure of the Banach $*$-algebra $A\widehat\otimes B$ requires further attention. We present some results in this direction in Section 3. We prove that the sum of two product ideals in $A\widehat{\otimes}B$ is closed and the same technique leads to a shorter proof of \cite[Theorem 3.8]{ass}. We further show that the minimal and maximal ideals in $A$ and $B$ generate their counterparts in $A\widehat\otimes B$. As a consequence, we obtain the lattice of closed ideals of $\B(H)\widehat{\otimes} \B(H)$.

%%%%%%%%%%%%%%%%%%%%%%%%%%%%%%%%%%%%%%%%%%%%%%%%%%%%%%%%%%%%%%%%%%%%%%%%%%%%%%%%
\section{Embedding Operator Space Projective Tensor Product into Second Duals}

  For operator spaces $V$ and $W$, a {\it jointly completely bounded bilinear map} (in short, j.c.b.) is a bilinear map $\phi:V \times W \rightarrow \mathbb{C}$ such that the maps $\phi_{n}:M_n(V)\times M_{n}(W)\rightarrow \M_{n^2}$ given by 
  $$
  \phi_n \big( (a_{ij}), (b_{kl}) \big )=\big(\phi(a_{ij},b_{kl})\big), \, \, n \in \mathbb{N}
  $$
  are uniformly bounded, and in this case we denote $\|\phi\|_{jcb}:= \text{sup}\{\| \phi_{n}\|:n \in \mathbb{N}\}$ \cite{bp}.    Also, a map $\phi:V \times W \rightarrow \mathbb{C}$ is said to be {\it completely bounded} (in short, c.b.) if the maps $\phi_{n}:M_n(V)\times M_{n}(W)\rightarrow \M_{n}$ given by 
  $$
  \phi_n \big( (a_{ij}), (b_{kl}) \big )=\big(\sum_{k=1}^{n}\phi(a_{ik},b_{kj})\big), \, \, n \in \mathbb{N}
  $$
are uniformly bounded, and then we write $\|\phi\|_{cb}:= \text{sup}\{\| \phi_{n}\|:n \in \mathbb{N}\}$. It is well known that $(V\oop W)^*$ and $(V\otimes_h W)^*$ are completely isometrically isomorphic to $\JCB(V\times W, \mathbb{C})$ and $\CB(V\times W,\mathbb C)$, respectively, where $\JCB(V\times W, \mathbb{C})$ (resp. $\CB(V\times W,\mathbb C)$) denotes the space of j.c.b.(resp. c.b.) bilinear maps \cite{bp,er1}. Every completely bounded map $\phi$ is jointly completely bounded with $\|\phi\|_{\text{jcb}} \leq \|\phi\|_{\text{cb}}$.

Recall that, for operator spaces $V$ and $W$, the {\it Haagerup norm} of an element $u \in M_n(V \otimes W)$, $n\in \N$, is defined by
$$
  \|u\|_h = \inf \{ \|x\| \|y\|: u= v \odot w, v \in M_{n,p}(V), w\in M_{p,n}(W), p\in \N \},  
$$
 where $v\odot w = (\sum_{k=1}^p v_{ik} \otimes w_{kj})_{ij}$.
%The Haagerup tensor product $V\otimes_h W$ is then the completion of $V\otimes W$ in the norm $\|\cdot\|_h$.
 The norms $\|\cdot\|_h,\|\cdot\|_{_{\wedge}}$ and $\|\cdot\|_\gamma$ on the tensor product $A\otimes B$ of two $C^*$-algebras $A$  and $B$ satisfy
$$
\|\cdot\|_h\, \leq\, \|\cdot\|_{_{\wedge}}\, \leq\, \|\cdot\|_\gamma.
$$

We first state an important result, whose proof can be found in \cite[\S 1.6.7]{blemer2}.

\begin{prop}\label{cb normal}
Let $V$ and $W$ be operator spaces, $E$ be a dual operator space, and $u:V\times W \ra E$ be a completely bounded bilinear map. Then $u$ admits a unique separately $w^*$-continuous extension $\tilde{u}: V^{**} \times W^{**} \ra E$, which is completely bounded with $\|u\|_{\text{cb}} = \|\tilde{u}\|_{\text{cb}}$.
\end{prop}

We now prove an embedding result for the Haagerup tensor product of operator spaces. It turns out that the operator space version is far much easier than the $C^*$-algebra case \cite[Theorem 4.1]{kumsi}, as observed below. It must be remarked that, if either $V$ or $W$ is finite dimensional, then $V^{**} \oh W^{**}$ is completely isometrically isomorphic to $(V \oh W)^{**}$ \cite[Corollary 9.4.8]{er1}.

\begin{thm}\label{haag-embed}
For  operator spaces $V$ and $W$, there is a canonical embedding of $V^{**} \oh W^{**}$ into $(V \oh W)^{**}$ which is a complete isometry.
\end{thm}
 
\begin{pf}
For the operator spaces $V^{**}$ and $W^{**},$ recall that $(V^{**} \oh W^{**})^*_\sigma$ denotes the subspace of $(V^{**}\oh W^{**})^*$ containing all the separately  $w^*$-continuous completely bounded bilinear forms on $V^{**} \times W^{**}$. By Proposition \ref{cb normal}, taking the map $u \ra \tilde{u}$, and $E$ as $\M_n$; one easily sees that $(V \oh W)^*$ is completely isometrically isomorphic to $ (V^{**} \oh W^{**})^*_\sigma.$ In particular, the normal Haagerup tensor product $V^{**} \otimes_{\sigma h} W^{**}$, which is defined as the operator space dual of  $(V^{**} \oh W^{**})^*_\sigma$, is completely isometrically isomorphic to $(V \oh W)^{**}$. Recall that, there also exists a completely isometric embedding \cite[\S 1.6.8]{blemer2} $ V^{**}\oh W^{**}  \hookrightarrow V^{**} \otimes_{\sigma h} W^{**}$. Hence, there is a completely isometric embedding  of $V^{**}\oh W^{**}$ into $ (V \oh W)^{**}$.
\end{pf}

\vspace*{2mm}
We now move on to analyze the embedding of biduals for operator space projective tensor product. This will need some preparations. Recall that an operator space $V$ is said to be {\it exact} if there exists a constant $K$ such that for any finite dimensional subspace $G\subset V$, there is an integer $n$, a subspace $\tilde{G} \subset \M_n$ and an isomorphism $u:G\ra \tilde{G}$ such that $\|u\|_{\text{cb}} \|u^{-1}\|_{\text{cb}} \leq K$. The smallest such constant is the exactness constant and is denoted by $ex(V)$. The matrix algebra $\M_n$ and in general, all nuclear $C^*$-algebras are simple examples of exact operator spaces with exactness constant 1. 
 
\begin{prop}\label{jcb normal}
Let $V$ and $W$ be exact operator spaces. Then every j.c.b. bilinear map $u:V\times W \rightarrow \mathbb C$ can be extended uniquely to a  separately $w^*$-continuous j.c.b. bilinear map $ \tilde{u}:V^{**}\times W^{**} \rightarrow \mathbb C$ such that $\|\tilde{u}\|_{jcb} \leq 2K \|u\|_{jcb}$, where $K=2\sqrt{2} \,ex(V)ex(W)$.  
\end{prop}

\begin{pf}
 Since $V$ and $W$ are both exact, by \cite[Theorem 0.4]{gpds} there exist bounded bilinear forms $u_1$ and $u_2$ on $V\times W$ such that $u=u_1 + u_2$ with $\|u_1\|_{\text{cb}} + \|u_2^t\|_{\text{cb}} \leq 2K\|u\|_{\text{jcb}}$, where $K=2\sqrt{2}\, ex(V)ex(W)$, and $u_2^t(w,v)= u_2(v,w)$. Using Proposition \ref{cb normal}, there exist unique separately $w^*$-continuous extensions $\tilde{u}_1: V^{**} \times W^{**} \ra \C$, and $\tilde{u}_2^t: W^{**} \times V^{**} \ra \C$ of $u_1$ and $u_2^t$,  which are completely bounded with $\|u_1\|_{\text{cb}} = \|\tilde{u}_1\|_{\text{cb}}$ and $\|u_2^t\|_{\text{cb}} = \|\tilde{u}_2^t\|_{\text{cb}}$. Note that $\tilde{u}_2^t$ is j.c.b, being c.b., so $\tilde{u}_2$ is also a j.c.b. bilinear form with $$\|\tilde{u}_2\|_{\text{jcb}} = \|\tilde{u}_2^t\| _{\text{jcb}} \leq \|\tilde{u}_2^t\|_{\text{cb}} = \|u_2^t\|_{\text{cb}}. $$ Also, it can be easily seen that $\tilde{u}_2$ is separately $w^*$-continuous, $\tilde{u}_2^t$ being separately $w^*$-continuous. Set $\tilde{u}= \tilde{u}_1 + \tilde{u}_2$. Then $\tilde{u}$ is a separately $w^*$-continuous j.c.b. bilinear form on $V^{**} \times W^{**}$ with $ \|\tilde{u} \|_{\text{jcb}} \leq 2K \|u \|_{\text{jcb}} $. Finally, some routine calculations show that $\tilde{u}$ is the unique extension of $u$.
\end{pf}

\vspace*{3mm}
In case of $C^*$-algebras $A$ and $B$, using the same techniques as that in Proposition \ref{jcb normal} and \cite[Lemma 3.1]{haagmusat}, one can easily prove that every j.c.b. bilinear map $u:A\times B \rightarrow \mathbb C$ can be extended uniquely to a  separately normal j.c.b. bilinear map $ \tilde{u}:A^{**}\times B^{**} \rightarrow \mathbb C$ such that $\|\tilde{u}\|_{jcb} \leq 2 \|u\|_{jcb}$. A priori, it is not clear why should this extension be norm preserving. However, we establish that, for $C^*$-algebras, via a completely different flavor, there exists a unique norm preserving separately normal extension. 

\begin{lem}{\label{lemma3}}
 Let $A$ and $B$ be von Neumann algebras and  $T: A \times B \rightarrow \mathbb{C}$ be a separately normal bilinear form. Then for each $n$, the map $T_n:M_n(A)\times M_n(B)\rightarrow \M_{n^2}$ defined as 
$$ T_n ( (a_{ij}),(b_{kl})) = (T(a_{ij},b_{kl})) $$
is separately normal.
\end{lem}

\begin{pf}
 For any $a=(a_{ij}) \in M_n(A)$ and for a fixed $b=(b_{ij}) \in M_n(B)$ we can write 
$$ T_n((a_{ij}), (b_{kl}))= \left(
\begin{array}{cccc}
T_{11}(a) & T_{12}(a) & \cdots & T_{1n}(a) \\
 T_{21}(a) & T_{22}(a) & \cdots & T_{2n}(a)\\
\vdots & \vdots & \vdots & \vdots\\
T_{n1}(a) & T_{n2}(a)  & \cdots & T_{nn}(a)
\end{array}\right), $$
  where $T_{kl}:M_n(A)\rightarrow \M_n$ maps $(a_{ij}) \rightarrow (T(a_{ij},b_{kl}))$. In order to show that the map $T_n$ is separately normal, we can equivalently show that the map $a\rightarrow \sum^{n}_{i,j=1} e_{ij} \otimes T_{ij}(a)$ from $M_n(A)$ into $ \M_n\otimes \M_n$ is normal. Let $(a_{\lambda})$ be an increasing net of positive elements in $M_n(A)$ such that $a_{\lambda} \stackrel{w^*}{\rightarrow}a$. Clearly, each $T_{ij}$ is normal, so $T_{ij}(a)$ is a weak limit of the net $(T_{ij}(a_{\lambda}))$ and thus $\sum e_{ij} \otimes T_{ij}(a)$ is a weak limit of $\sum e_{ij} \otimes T_{ij}(a_{\lambda})$. Hence the result.
\end{pf}

\begin{prop}{\label{jcbnormal2}}
Let $A$ and $B$ be $C^\ast$-algebras and $\phi:A\times B \rightarrow \mathbb C$ be a j.c.b. bilinear form. Then $\phi$ admits a unique separately normal j.c.b. bilinear extension $ \tilde{\phi}:A^{**}\times B^{**} \rightarrow \mathbb C$ such that $||\tilde{\phi}||_{jcb} = ||\phi||_{jcb}$. 
\end{prop}

\begin{pf}
 Since $\phi:A\times B \rightarrow \mathbb C$ is a continuous bilinear form, there exists a unique separately normal bilinear form $ \tilde{\phi}:A^{**}\times B^{**} \rightarrow \mathbb C$ with $\| \tilde{\phi}\| = \| \phi \|$ \cite[Corollary 2.4]{haag}. Let $n$ be a positive integer. Consider the map $\tilde{\phi}_n:M_n(A^{**})\times M_n(B^{**})\rightarrow \M_{n^2}$ defined as 
$$
\tilde{\phi}_n ( (a_{ij}),(b_{kl})) = (\tilde{\phi}(a_{ij},b_{kl})).
$$
We claim that $\| \tilde{\phi}_n\| \leq \| \phi\|_{jcb}$.  Consider any $a \in M_n(A^{**}),\, b \in M_n(B^{**})$ with $\|a\| = \|b\| = 1$. Using the fact that the unit ball of $M_n(A)$ is $w^*$-dense in the unit ball of $M_n(A^{**})$, we obtain a net $(a_{\lambda}) \in M_n(A) $ which is $w^*$-convergent to $a$ with $\|a_{\lambda}\| \leq 1$ and a net $(b_{\mu}) \in M_n(B) $ which is $w^*$-convergent to $b$ with $\|b_{\mu}\| \leq 1$. By Lemma \ref{lemma3}, $\tilde{\phi}_n$ is separately normal, so for a fixed $k \in \mathbb{N}$, 
$$
 \| \tilde{\phi}_n(a,b)\|  \leq (1+1/k)^2 \| \tilde{\phi}_n(a_{\lambda},b_{\mu})\|,$$
 for some $\lambda$ and $\mu$. This further gives
$$ \| \tilde{\phi}_n(a,b)\| \leq (1+1/k)^2 \| \phi_n\|, \quad \forall\, k \in \mathbb{N}.$$
Thus  $\| \tilde{\phi}_n\| \leq \|\phi_n\| \leq \| \phi\|_{jcb}$, and this is true for all $n\in \N$,  giving that $\tilde{\phi}$ is j.c.b. with $\|\tilde{\phi}\|_{\text{jcb}} \leq \| \phi \|_{\text{jcb}}$. Also $\| \phi \|_{\text{jcb}} \leq \| \tilde{\phi} \|_{\text{jcb}}$, $\phi$ being the restriction of $\tilde{\phi}$, and hence the result.
\end{pf}

\vspace*{3mm}
We next prove a result which is an operator space version of \cite[Lemma 5.3]{kumsi}, and whose proof is largely inspired by the same.
\begin{lem}\label{approx}
 Let $V$ and $W$ be operator spaces with $V\subset \B(H)$ and $W \subset \B(K)$. Then the unit ball of $\CB_\sigma(V\times W,\C)$ is $w^*$-dense in the unit ball of $\CB(V \times W, \C)$, where  $\CB_\sigma(V\times W,\C)$ denotes the space of all separately $w^*$-continuous c.b. bilinear forms on $V\times W$.
\end{lem}

\begin{pf}
Let $B_1$ and $B_2$ denote the unit balls of $\CB_\sigma(V\times W,\C)$ and $\CB(V \times W, \C)$, respectively. Let if possible, there exist a $\phi$ in $B_2$ such that $\phi \notin \bar{B_1}^{w^*}$, where $\bar{B_1}^{w^*}$ denotes the $w^*$-closure of $B_1$ in $B_2$. Using a consequence of  Hahn Banach separation Theorem \cite[Theorem 3.7]{rudin}, we obtain a $w^*$-continuous linear functional $\Phi:\CB(V \times W, \C) \ra \C$ such that $|\Phi (\psi)| \leq 1$  for all $ \psi \in B_1$ and $\Phi(\phi) >1$. Now $\Phi$ can be identified with an element $u$ of $V\oh W$, being a $w^*$-continuous functional on $(V \oh W)^*$. Therefore, there is a $u \in V \oh W$, $\|u\|_h >1 $, such that $| \psi(u)| \leq 1$ for all $\psi \in B_1$, and $\phi(u) >1$.

It is well known that there is an isometric embedding of $\B(H) \oh \B(K) $ into $\CB (\B(K,H))$ \cite[Theorem 4.3]{smith}. Using the injectivity of the Haagerup tensor product, we get an isometric embedding (need not be algebraic), say $\theta$, of $V\oh W$  into $\CB (\B(K,H))$ given by $\theta (v\otimes w)(T)=vTw$. Since $\|\theta(u)\|_{\text{cb}} > 1$, for some $n\in \mathbb{N}$,  $\| (\theta(u))_n\| > 1$. So, there exists $(x_{ij}) \in M_n(\B(K,H))$ with $\|(x_{ij})\| =1$ such that  
\[\| (\theta(u)x_{ij})\|=\| (\theta(u))_n (x_{ij})\| >1.\]
Now, we can choose unit vectors $\xi \in K^n$ and $\eta \in H^n$ such that
\[ | \langle (\theta(u)x_{ij})\xi, \eta \rangle | > 1. \]
Define $\psi: V \times W \ra \C$ as 
\[  \psi(v, w) = \langle (\theta(v\otimes w)x_{ij})\xi, \eta \rangle.\]
Clearly $\psi$ belongs to $B_1$, which together with the relation
$$|\psi(u)| = |\langle (\theta(u)x_{ij})\xi, \eta \rangle| >1,$$
give  a contradiction. Hence $B_1$ is $w^*$-dense in $B_2$.
\end{pf}

\vspace*{3mm}
Now, with all the necessary ingredients at our disposal, we are ready to prove the main result of this section; which is an operator space analogue of \cite[Theorem 5.1]{kumsi}, and the proof presented here borrows ideas from the same. Let us first define the required embedding. Using Proposition \ref{jcb normal}, we have a map $\chi: (V \oop W)^* \ra (V^{**}\oop W^{**})^* $ with $\| \chi\| \leq 2K$. Define $\mu :=\chi^{*} \circ i: V^{**}\widehat{\otimes}W^{**} \ra
(V\widehat{\otimes}W)^{**}$, where $i:V^{**}\widehat{\otimes}W^{**}\ra(V^{**}\widehat{\otimes}W^{**})^{**}$ is the canonical completely isometric embedding. Then, we have the following:

\begin{thm}\label{opembed}
  For exact operator spaces $V$ and $W$, the embedding $\mu$ of $V^{**}\widehat{\otimes}W^{**}$ into
  $(V\widehat{\otimes}W)^{**}$ satisfies 
  $$
   \frac{1}{2K} \|u\| \leq \|\mu(u)\| \leq 2K\|u\|, \, \forall \, u \in V^{**}\widehat{\otimes}W^{**},
  $$
 where $K=2\sqrt{2} \,ex(V)ex(W)$. In particular, $\mu$ has a continuous inverse.
\end{thm}

\begin{pf}
By the definition of $\mu $, the inequality on the R.H.S. is obvious. For the other inequality, consider any $u\in V^{**}\otimes W^{**} $  with $\|u\|_\wedge=1 $. By Hahn Banach Theorem, there exists a j.c.b. bilinear map $\phi: V^{**}\times W^{**}  \ra \mathbb{C}$ such that $\|\phi\|_{jcb}=1$ and $\phi(u)=1$.  By \cite[Theorem 0.4]{gpds}, $\phi$ can be decomposed as $\phi =\phi_1 +\phi_2$, where $\phi_1$ and $\phi_2$ are bounded bilinear forms with $\|\phi_1\|_{\text{cb}} \leq K$ and $\|\phi^t_2\|_{\text{cb}} \leq K$, where $\phi^t_2(b,a)=\phi_2(a,b)$. Now consider any $\epsilon >0$. For the $w^*$-open sets $\{\theta \in B_1(\CB(V^{**}\times W^{**},\C)): |(\theta -\phi_1/K)(u)|< \epsilon \}$ and  $\{\zeta \in B_1(\CB(W^{**}\times V^{**},\C)): |(\zeta -\phi_2^t/K)(u^t)|< \epsilon \}$, using Lemma \ref{approx}, we get $\Phi_1 \in \CB_\sigma(V^{**}\times W^{**},\C)$ and  $\Phi_2 \in \CB_\sigma(W^{**}\times V^{**},\C) $ with $\|\Phi_j\|_{\text{cb}} \leq 1$, $j=1,2$,  such that
%The c.b. bilinear forms $\phi_1/K$ and $\phi^t_2/K$ may be weak*-approximated on the  elements $u$ and $u^t$ by normal c.b. bilinear forms 
  \[ |\phi_1(u)-K\Phi_1(u)| < K\epsilon , \quad |\phi^t_2(u^t)-K\Phi_2(u^t)| < K\epsilon\]
which further give
\begin{equation} \label{i}
 |\phi_1(u)-K\Phi_1(u)| < K\epsilon , \quad |\phi_2(u)-K\Phi^t_2(u)| < K\epsilon. 
\end{equation}
  Now $\Phi_1$ and $\Phi_2$ both are j.c.b., being c.b.. Also $\Phi^t_2$ is a separately $w^*$-continuous j.c.b. form on $V^{**} \times W^{**}$ with  $\|\Phi^t_2\|_{\text{jcb}}= \|\Phi_2\|_{\text{jcb}}$, so $\Phi =\Phi_1 +\Phi^t_2$ is a j.c.b. map. Let $\psi_1, \, \psi_2$ be the restrictions of $\Phi_1$ and $\Phi^t_2$ to $ V \times W$, then these are j.c.b. bilinear 
maps. Thus, by the definition of $\chi,\, \Phi_1 = \chi(\psi_1), \, \Phi^t_2 = \chi(\psi_2)$. Set $\psi =\psi_1 + \psi_2$. Then $\psi$ is a j.c.b. bilinear map and thus it is a continuous  linear functional on $V\widehat{\otimes} W$ with $\| \psi\| \leq 2$. Further,
$$
\mu(u)(\psi)= \chi^* i(u) (\psi)=i(u)(\chi\psi) = (\chi \psi)(u) = \Phi(u)
$$
which, along with (\ref{i}), give $\| \mu(u) \| \geq 1/2K$.
\end{pf}
 
%\vspace*{3mm}
\begin{rem} If the extension of \cite[Conjecture $0.2^\prime$]{gpds} is true for $\M_n$-valued bilinear functions, which is not known to us, then we can prove that $\mu$ is completely bounded. Indeed, if the conjecture were true, then using the same argument as that in the proof of Proposition \ref{jcb normal}, one can prove that for exact operator spaces $V$ and $W$, every j.c.b. bilinear map $u:V\times W \rightarrow \M_n$ can be extended uniquely to a separately $w^*$-continuous j.c.b. bilinear map $ \tilde{u}:V^{**}\times W^{**} \rightarrow \M_n$ such that $\|\tilde{u}\|_{jcb} \leq 2K \|u\|_{jcb}$, for some constant $K$ independent of $n$. Now, to show that $\mu$ is completely bounded, it is sufficient to show that $\chi$ is so, and in that case $\|\mu\|_{\text{cb}} \leq \| \chi\|_{\text{cb}}$. Note that, for $n\in \mathbb{N}$, we have the following commutative diagram:
\[ \xymatrix{ M_n(\JCB(V\times W,\C)) \ar[r]^{\chi_n} \ar[d]_{i} & M_n(\JCB(V^{**}\times W^{**},\C)) \ar[d]^{i^\prime} \\
      \JCB(V\times W,\M_n) \ar[r]^{\chi^\prime} & \JCB(V^{**}\times W^{**},\M_n)}\]
In the above diagram, $\|\chi^\prime\| \leq 2K$, $K$ being independent of $n$, and both $i$ and $i^\prime$ are complete isometric isomorphisms. So, $\|\chi_n\| \leq 2K$, and this is true for all $n\in \N$. Thus $\chi$ is completely bounded with $\|\chi\|_{\text{cb}} \leq 2K$.
\end{rem}

For $C^*$-algebras $A$ and $B$, using Proposition \ref{jcbnormal2} and the techniques of Theorem \ref{opembed}, one can prove the following: 

\begin{thm}\label{embed}
For $C^*$-algebras $A$ and $B$, there is a  canonical bi-continuous embedding $\mu$ of $A^{**}\widehat{\otimes}B^{**}$ into $(A\widehat{\otimes}B)^{**}$ which  satisfies 
  $$   \frac{1}{2} \|u\| \leq \|\mu(u)\| \leq \|u\|, \quad \forall \, u \in A^{**}\widehat{\otimes}B^{**}. $$
\end{thm}

As an application of the above result, we prove an equivalence between the Haagerup norm and the operator space projective norm for tensor product of $C^*$-algebras. This result has already been proved by Kumar and Sincalir \cite[Theorem 7.4]{kumsi}. However we use a different and rather simple technique to prove the same. We first need the following easy result dealing with the injectivity of the projective norm.

\begin{lem}\label{fd}
If $A_0$ and $B_0$ are both finite dimensional $C^*$-subalgebras of the $C^*$-algebras $A$ and $B$, then $A_0 \oop B_0$ is a closed $*$-subalgebra of $A\oop B$.
\end{lem}

\begin{pf}
 Since $A_0$ and $B_0$ are both finite dimensional $C^*$-subalgebras of $A$ and $B$, there are conditional expectations $P_1$ and $P_2$ from $A$ and $B$ onto $A_0$ and $B_0$, respectively, with $\|P_1\|_{\text{cb}} = \|P_2\|_{\text{cb}} =1$, see \cite[II.6.10.4]{blackadar}.  Then $P_1 \oop P_2$ is a projection of $A \oop B$ onto $A_0 \oop B_0$ with $\|P_1 \oop P_2\| \leq 1$. For the inclusion map $i: A_0 \oop B_0 \ra A\oop B$, the composition $(P_1 \oop P_2)\circ i$ agrees with the identity map on $A_0 \otimes B_0$, so that for any element $x \in A_0 \otimes B_0$,
\begin{eqnarray*}
 \|x\|_{A_0\oop B_0} & = & \|((P_1 \oop P_2)\circ i )(x)\|_{A_0\oop B_0} \\
                     & \leq & \|i(x)\|_{A\oop B} \\
                      & \leq & \|x\|_{A_0\oop B_0}.
\end{eqnarray*}
Hence $i$ is an isometry, giving $A_0\oop B_0$ as a closed subalgebra of $A\oop B$.
\end{pf}

\vspace*{3mm}
Recall that, a $C^*$-algebra $A$ is said to be {\it $n$-subhomogenous} if each irreducible representation of $A$ has dimension less than or equal to $n$, and {\it subhomogenous} if it is $n$-subhomogenous for some $n \in \mathbb{N}$. It is known that a $C^*$-algebra $A$ is $n$-subhomogenous if and only if $A^{**}$ does not contain a $C^*$-subalgebra  isomorphic to $\M_{n+1}$.

\begin{thm}
 For $C^*$-algebras $A$ and $B$, the Haagerup norm $\|\cdot\|_h$ is equivalent to the operator space projective tensor norm $\|\cdot \|_\wedge$ on $A\otimes B$ if and only if $A$ and $B$ are both subhomogenous.
\end{thm}

\begin{pf}
 Let us assume that $A$ and $B$ are both infinite dimensional, and
\[ \|x\|_\wedge \leq c\|x\|_h,\quad  \forall \, x \in A \otimes B, \]
for some constant $c$, that is,  the canonical map $j:A\oh B \ra A\oop B$ is continuous with $\|j\|\leq c$. We first claim that
$$\|x^{**}\|_\wedge \leq 2c\|x^{**}\|_h,\qquad \forall x^{**} \in A^{**} \otimes B^{**}.$$
In other words,  the identity map $J:A^{**}\otimes B^{**} \ra A^{**}\oop B^{**}$ is continuous with respect to the Haagerup norm.  Using Theorem \ref{embed} and Theorem \ref{haag-embed}, we have a bi-continuous canonical embedding $\mu:A^{**}\oop B^{**} \ra (A \oop B)^{**}$, and   a canonical completely isometric embedding $\zeta:A^{**}\oh B^{**} \ra (A \oh B)^{**}$, respectively. For any $a^{**} \otimes b^{**} \in A^{**}\otimes B^{**}$ and $f\in  (A\oop B)^*$, we have
$$ \mu J (a^{**} \otimes b^{**})(f)  = (\chi f)(a^{**} \otimes b^{**}) = \tilde{f}(a^{**} \otimes b^{**}),$$
where $\tilde{f}:A^{**}\times B^{**} \ra \C$ is the unique separately normal j.c.b. extension of $f:A\times B \ra \C$. Also 
\[ j^{**} \zeta (a^{**} \otimes b^{**})(f) = \zeta (a^{**} \otimes b^{**})(j^*(f)) =\widetilde{j^*(f)} (a^{**} \otimes b^{**}),\]
where $\widetilde{j^*(f)}:A^{**}\times B^{**} \ra \C$ is the unique separately normal c.b. extension of $j^*(f):A\times B \ra \C$. Note that, $\widetilde{j^*(f)}$ is also a j.c.b. extension of $f$, so by uniqueness,  $\widetilde{j^*(f)} = \widetilde{f}$, which gives $\mu J = j^{**}\zeta$ on $A^{**}\otimes B^{**}$. Using the bi-continuity of $\mu$, we get $J=\mu^{-1}j^{**}\zeta$, with
\begin{equation}\label{III}
\|J\| \leq \|\mu^{-1}\| \|j^{**}\| \|\zeta\| \leq 2c,	
\end{equation}
which proves our first claim.

Let $A^{**}$ contains an isomorphic copy (not necessarily unital) of $\mathbb{M}_n$, for some $n\in \mathbb{N}$. Then $B^{**}$ also contains a copy of $l_n^{\infty}$. Using the injectivity of Haagerup norm and Lemma \ref{fd}, $\mathbb{M}_n\oh l_n^{\infty}$ and $\mathbb{M}_n\oop l_n^{\infty}$ embed isometrically in $A^{**} \oh B^{**}$ and $A^{**} \oop B^{**}$ respectively. Let $\{e_{ij}\}$ denote the standard matrix units, then using \cite[Lemma 3.1]{kumsi}, and (\ref{III}), we have
\begin{eqnarray*}
 n^{1/2} = \| \Sigma^n_{j=1} \,e_{1j} \otimes e_{jj}\|_h & \leq &  \| \Sigma_j\, e_{1j} \otimes e_{jj}\|_\wedge\\
  & =& \| \Sigma_j\, e_{j1} \otimes e_{jj}\|_\wedge \\
  & \leq & 2c \| \Sigma_j\, e_{j1} \otimes e_{jj}\|_h \\
  & = & 2c.    
\end{eqnarray*}
So $A^{**}$ cann't contain an isomorphic copy of $\mathbb{M}_n$ for $n> 4c^2$ which shows that $A$ is $4c^2$-subhomogenous. A similar argument gives that $B$ is also $4c^2$-subhomogenous. 

The other implication is easy. 
\end{pf}

%%%%%%%%%%%%%%%%%%%%%%%%%%%%%%%%%%%%%%%%%%%%%%%%%%%%%%%%%%%%%%%%%%%%%%%%%%%%%%%%%%%%%%%
\section{\texorpdfstring{Ideal Structure of $A\widehat{\otimes}B$}{Ideal Structure of A oop B}}

The operator space projective tensor norm is symmetric, associative and projective but not injective \cite{er1}. For $C^*$-algebras $A$ and $B$, $A\widehat{\otimes}B$ is a Banach $*$-algebra with the natural isometric involution given by $*:a\otimes b\rightarrow a^*\otimes b^*$ \cite{kumar}. This property is in contrast to the Haagerup norm, where  the natural involution on $A\otimes_h B$ is an isometry if and only if $A$ and $B$ are commutative \cite{kumar2}. This section is devoted to a systematic study of the ideal structure of this Banach $*$-algebra. If $K$ and $L$ are closed ideals of $A$ and $B$, where $A$ and $B$ are $C^\ast$-algebras, then $K\widehat{\otimes}L$ is a closed $*$-ideal of $A\widehat{\otimes}B$ \cite[Theorem 5]{kumar}, which is termed as {\it product ideal}. Allen, Sincalir and Smith \cite{ass} proved that sum of two product ideals in the Haagerup tensor product is again a closed ideal. In this section we discuss its analogue for the operator space projective tensor product, whose techniques also give a shorter proof of \cite[Theorem 3.8]{ass}. It must be mentioned that all the results in this section hold true for  closed $*$-ideals also. Through out this section $A$ and $B$ denote the $C^*$-algebras, unless otherwise stated. We first state an elementary result, a proof of which for the Banach space projective norm can be found in \cite{kan}.

\begin{lem}\label{app}
 If Banach algebras $A$ and $B$ both possess bounded approximation identities, then for any subcross norm $\alpha$, $A\otimes^{\alpha}B$ possesses a bounded approximation identity, where  $A\otimes^{\alpha}B$ is the completion of the algebraic tensor product $A\otimes B$ with respect to $\alpha$ norm.
\end{lem}

\begin{prop}\label{prop3.1}
    Let $I_1, I_2$ and $J_1, J_2$ be the closed ideals of $A$ and $B$ respectively. Then $I_1 \widehat\otimes J_1 + I_2 \widehat\otimes J_2 $ is a closed $*$-ideal of $A \widehat\otimes B$. 
\end{prop}

\begin{pf}
 By Lemma \ref{app} and \cite[Theorem 5]{kumar}, it follows that $I_1 \widehat\otimes J_1$  and $ I_2 \widehat\otimes J_2 $ are closed $*$-ideals, both having bounded approximation identities. Using the fact that sum of two closed ideals is closed if any one of them has bounded approximate identity \cite[Prop 2.4]{dixo}, we obtain the required result.
\end{pf}

\begin{rem}\label{haagideal}
The above proposition is true for Haagerup norm and  Banach space projective norm as well. In particular, this gives a shorter proof of \cite[Theorem 3.8]{ass}. 

\end{rem}

As a direct consequence of Proposition \ref{prop3.1}, we next show that the  operator space projective tensor product is distributive over the finite sums of closed ideals. 

\begin{cor}\label{prop1}
If ${M_i}$ and $N_i$, $i=1,2,\dots,n$, are closed ideals in $A$ and $B$ respectively, then
\begin{enumerate}
\item $A \widehat{\otimes}(\Sigma^{n}_{i=1} N_i) =\sum^{n}_{i=1}(A \widehat{\otimes} N_i)$,
\item $(\Sigma^{n}_{i=1} M_i)\widehat{\otimes} B = \sum^{n}_{i=1}(M_i \widehat{\otimes} B)$.
\end{enumerate}
\end{cor}

\begin{pf}
We shall only prove  the first part, and the proof for (2) follows on the same lines. Using \cite[Theorem 5]{kumar}, each $A\oop N_i$ is a closed ideal of $A\oop (\Sigma_i N_i)$, so it is easy to see that 
\[
 A \oop (\Sigma_i N_i) \supseteq \Sigma_i(A \oop N_i).
\]
For the other containment, note that $A \otimes (\Sigma_i N_i) \subseteq  \Sigma_i (A\oop N_i)$. By \cite[Theorem 5]{kumar} and Proposition \ref{prop3.1},  $A \oop (\Sigma_i N_i)$ and $\Sigma_i (A \oop N_i)$ are both closed in $A\oop B$. So $A \oop (\Sigma_i N_i) \subseteq  \Sigma_i A\oop N_i$, proving the result.
\end{pf}

\vspace*{2mm}

We would like to remark that Allen, Sinclair and Smith,  proved the analogue of the above result for Haagerup tensor product \cite[Proposition 2.9]{ass}. However, their method was more technical. Using Remark \ref{haagideal} and the same argument as in the proof of the above result, a much shorter and simpler proof can be provided for the same. 

In case of $C^\ast$-algebras, again using Proposition \ref{prop3.1}, we have the following modified version of \cite [Proposition 7.1.7]{er1}: 

\begin{prop}\label{prop3.2}
Let $A, A_1, B$ and $B_1$ be $C^\ast$-algebras. Given the (complete) quotient mappings $\phi: A \rightarrow A_1$ and $\psi: B \rightarrow B_1$, the corresponding mapping $ \phi \otimes \psi : A\otimes B\rightarrow A_1\otimes B_1$ extends to a (complete) quotient mapping $ \phi \widehat\otimes \psi:A  \widehat\otimes B \rightarrow A_1 \widehat\otimes B_1 $.
Further, \[\ker (\phi \,\widehat\otimes  \psi) = \ker \phi\, \widehat\otimes \,B + A \,\widehat\otimes \ker \psi.\]
\end{prop}

\begin{pf}
   From \cite [Proposition 7.1.7]{er1}, we know that \[\ker(\phi \widehat\otimes  \psi) = (\ker \phi \otimes B + A \otimes \ker \psi)^-,\] so it enough to check that \[ (\ker \phi \otimes B + A \otimes \ker \psi)^- = \ker \phi \,\widehat\otimes\, B + A \,\widehat\otimes \ker \psi.\] 
   Note that  $\ker \phi \, \widehat\otimes B $ and $A \,\widehat\otimes \ker \psi$ are closed ideals of $A\widehat\otimes B$ \cite[Theorem 5]{kumar} and they can be realized as the closure of $\ker \phi \otimes B$ and $A \otimes \ker \psi $ in $A\widehat\otimes B$. The result now follows easily using the fact that $\ker \phi \,\widehat\otimes B + A \,\widehat\otimes \ker \psi$ is closed. 
\end{pf}

\vspace*{2mm}
In \cite{jain}, we proved that the canonical map $i:A\widehat{\otimes}B \ra A\otimes_{\text{min}}B$ is injective. Making repeated use of this result along with some techniques of Allen, Sinclair and Smith \cite{ass}, we will now study the ideal structure of $A\widehat{\otimes}B$ in terms of the ideal structures of $A$ and $B$.

\begin{prop}\label{prop2}
   Let $I$ be a non-zero closed ideal of $A\widehat{\otimes}B$. Then $I$ contains a non-zero elementary tensor and a non-zero product ideal.
\end{prop}

\begin{pf}
  Let $I_{\min}$ denote the min-closure of $I$ in $A\otimes_{\text{min}}B$, i.e., $I_{\min}$ is the closure of $i(J)$ in $A\omin B$. Then $I_{\text{min}}$ is a non-zero ideal of $A{\otimes}_{\text{min}}B$ \cite[Corollary 1]{jain}, and thus contains a non-zero elementary tensor \cite[Proposition 4.5]{ass}, say $a\otimes b$, which also lies in $I$ \cite[Theorem 6]{kumar}. Let $K$ and $L$ be the non-zero closed ideals in $A$ and $B$ generated by $a$ and $b$. Then clearly $I$ contains the product ideal $K\widehat{\otimes}L$.
\end{pf}

\begin{thm}\label{theo1}
   The Banach $*$-algebra $A\widehat{\otimes}B$ is simple if and only if $A$ and $B$ are simple.
\end{thm}

\begin{pf}
  Let $I$ be a non-zero closed ideal of $A\widehat{\otimes}B$. Then by Proposition \ref{prop2}, $I$ contains a non-zero product ideal  $K\widehat{\otimes}L$, where $K$ and $L$ are non-zero ideals of $A$ and $B$ respectively. But $A$ and $B$ are simple so $K=A$ and $L=B$.  Thus $A\widehat{\otimes}B$ is simple.

 For the reverse implication, let if possible $A$ be not simple. Then it contains a non-trivial closed ideal, say $I$, which gives rise to a non-zero closed ideal $I\oop B$ of $A\oop B$. Now Proposition \ref{prop3.2} gives  an isomorphism between the spaces $(A\oop B) / (I\oop B)$ and $(A/I) \oop B$, which implies $I\oop B$ is proper in $A\oop B$. Thus $I\oop B$ is a non-trivial closed ideal of $A\oop B$, which contradicts the fact that $A\oop B$ is simple. Similarly one can prove that $B$ is simple.
\end{pf}

\begin{thm}\label{theo2}
  Let $A$ and $B$ be $C^*$-algebras with $A$ as simple, then every closed ideal in $A\widehat{\otimes}B$ has the form $A\widehat{\otimes}L$ for some closed ideal $L$ in $B$.
\end{thm}

\begin{pf}
%  Let $I$ be a non-zero closed ideal in $A\widehat{\otimes}B$. Then by Proposition \ref{prop2}, $I$ contains a non-zero product ideal of the form $A\widehat{\otimes}L_1$, $L_1$ being a non-zero closed ideal of $B$. If $L_1$ and $L_2$ are closed ideals of $B$ with $A\widehat{\otimes}L_1$ and $A\widehat{\otimes}L_2$ contained in $I$, then by Corollary \ref{prop1}, 
 % \[ A\widehat{\otimes}(L_1 + L_2)=A\widehat{\otimes}L_1 + A\widehat{\otimes}L_2 \subseteq I. \]
  % Let $\mathcal F$ be the family of closed ideals $L$ of $B$ such that $A\widehat{\otimes}L\subseteq I$. For any chain $\mathcal P=\{L_i:i\in\Lambda\}$ in $\mathcal F$, $J=\{\sum_{\text{finite}} x_i:x_j\in L_j,j\in\Lambda\}$ is an upper bound of $\mathcal P$. Thus by Zorn's lemthere is a largest closed ideal $L \subseteq B$ such that $A\widehat{\otimes}L \subseteq I$.
  Let $K$ be a non-zero closed ideal in $A\widehat{\otimes}B$. By Proposition \ref{prop2}, since $A$ is simple, $K$ contains a non-zero product ideal of the form $A\widehat{\otimes}L_1$, $L_1$ being a non-zero closed ideal of $B$. 
   Consider the non-empty family $\mathcal F$ of closed ideals $L$ of $B$ such that $A\widehat{\otimes}L\subseteq K$. Let $\mathcal P=\{L_i:i\in\Lambda\}$ be a chain in $\mathcal F$. Note that, by Corollary \ref{prop1}, $\mathcal{F}$ is closed under finite sums. So, $J=\{\sum_{\text{finite}} x_i:x_j\in L_j,j\in\Lambda\}$ is an upper bound of $\mathcal P$ in $\mathcal{F}$. Thus by Zorn's Lemma, there is a largest closed ideal $L \subseteq B$ such that $A\widehat{\otimes}L \subseteq K$.

   Consider, the quotient map $1\otimes \pi:A\widehat{\otimes}B \rightarrow A\widehat{\otimes}(B/L)$ with kernel $A\widehat{\otimes}L$. Then $\tilde{K}= (1\otimes \pi)(K)$ is a closed ideal of $A\widehat{\otimes}(B/L)$. It is sufficient to show that $\tilde{K}$ is a zero ideal, as in that case $K \subseteq \text{ker}(1 \otimes \pi) = A\widehat{\otimes}L$. If $\tilde{K}$ is non-zero, then it would contain a non-zero elementary tensor say $a \otimes(b+L)= (1 \otimes \pi)(a\otimes b)$, where $a \otimes b \in K$. Let $N$ be the closed ideal in $B$ generated by $b$. Since $A$ is simple, $K$ contains the closed ideal $A\widehat{\otimes}N$. But $A\widehat{\otimes}N$ is not contained in $A\widehat{\otimes}L$, which contradicts the maximality of $L$. Thus $\tilde{K}$ is zero ideal and hence the result.
\end{pf}

\begin{prop}\label{min}
A closed ideal $J$ in $A\widehat{\otimes}B$ is minimal if and only if there exist minimal closed ideals $K \subseteq A$ and $L  \subseteq B$ such that $J = K\widehat{\otimes}L$.
\end{prop}

\begin{pf}
  Let $J$ be minimal in $A\widehat{\otimes}B$. By Proposition  \ref{prop2}, there is a non-zero product ideal $K\widehat{\otimes}L$ contained in $J$. Since $J$ is minimal, $J =  K\widehat{\otimes}L$, and it is clear that $K$ and $L$ must be  minimal in $A$ and $B$ respectively.
 
  Conversely, let $K$ and $L$ be minimal closed ideals. Then they both are simple $C^*$-algebras. By Theorem \ref{theo1}, $K\widehat{\otimes}L$ is simple and thus contains no proper non-zero closed ideal of $A\widehat{\otimes}B$. Hence it is minimal.
\end{pf}

\begin{thm}\label{max}
A closed ideal $J$ is maximal in $A\widehat{\otimes}B$ if and only if there exist maximal closed  ideals $M$ in $A$ and $N$ in $B$ such that \[J = A \widehat{\otimes} N + M \widehat{\otimes} B.\] 
\end{thm}

\begin{pf}
  Let $M$ and $N$ be  maximal ideals of $A$ and $B$ respectively. Note that, by Proposition \ref{prop3.1}, $J = A\widehat{\otimes}N + M\widehat{\otimes}B$ is a closed ideal of $A\widehat{\otimes}B$. Also if $\pi_1 : A  \rightarrow A/M$ and $\pi_2 : B \rightarrow B/N$ are quotient maps, then by Proposition \ref{prop3.2} $ J$ is equal to $\text{ker}(\pi_1 \otimes \pi_2)$, and there is an isomorphism between $(A\widehat{\otimes}B)/J$ and $(A/M)\widehat {\otimes} (B/N)$. By Theorem \ref{theo1}, $(A\widehat{\otimes}B)/J$ is a simple Banach*-algebra. Thus J is maximal in $A\widehat{\otimes}B$.

Conversely, let $J$ be a maximal ideal of $A\widehat{\otimes}B$. Let $J_{\text{min}}$ be the min-closure of $J$ in $A\otimes_{\text{min}}B$. Then $J_{\text{min}}$ is a non-zero closed ideal of $A\otimes_{\text{min}}B$ \cite{jain} and it is proper since $J_{\text{min}}=A\otimes_{\text{min}}B$ would imply $J=A\widehat{\otimes}B$ \cite{kumar}. Let $ \pi : A\otimes_{\text{min}}B \rightarrow B(H)$ be an irreducible representation annihilating $J_{\text{min}}$. Since the canonical map $i:A\widehat{\otimes}B \rightarrow A\otimes_{\text{min}} B$ is bounded $*$-homomorphism, we get a $*$-representation $\tilde{\pi}=\pi \circ i$ of $A\widehat{\otimes} B$ on $H$ such that $\tilde{\pi} (J) = \{0\}$. By \cite[Lemma IV.4.1]{take},  there exist commuting representations $\pi_1$ and $\pi_2$ of $A$ and $B$ on $H$, respectively such that $ \tilde{\pi} (a\otimes b) = \pi_1 (a) \pi_2 (b),\,\, \forall \,a\in A,\,b\in B $. Let $M= \text{ker}\, \pi_1,\,N= \text{ker} \,\pi_2$ and $I=A\widehat{\otimes}N + M\widehat{\otimes}B$. Clearly $\tilde{\pi}(M\widehat{\otimes} B)=\{0\}= \tilde{\pi}(A\widehat{\otimes} N)$, which gives $\tilde{\pi}(J+I)=\{0\}$. So $J+I$ is a proper ideal of $A\widehat{\otimes}B$, which  by maximality of $J$, gives $I \subseteq J$. For the reverse inclusion, using  Proposition \ref{prop3.2}, there is a quotient map $ q : A \widehat{\otimes}B\rightarrow (A/M) \widehat{\otimes} (B/N)$  with kernel $I$. It is sufficient to show that $q(J)=\{0\}$. Now, the representations $\pi_1$ and $\pi_2$ induce faithful commuting representations $\tilde{\pi_1}$ of $A/M$ and $\tilde{\pi_2}$ of $B/N$ on $H$. Using \cite[Proposition IV.4.7]{take} and the fact that the canonical map $i:A\widehat{\otimes}B \rightarrow A\otimes_{\text{max}} B$ is a bounded $*$-homomorphism, there exists a representation $\pi_0$ of $(A/M) \widehat{\otimes}(B/N)$ on $H$ such that $\pi_0(x\otimes y) = \tilde{\pi}_1 (x)\tilde{\pi}_2(y),\, \forall \, x\in A/M, y\in B/N $. It is easy to verify that, $ \tilde{\pi}$ and $\pi_0 \circ q $ agree on $A\otimes B$, which by continuity gives $\pi_0 (q(J))=0 $. Now, $\pi$ is an irreducible representation, so $\tilde{\pi_1}$ and $\tilde{\pi_2}$ are both faithful factor representations with commuting ranges. Using \cite[Proposition IV.4.20]{take}, $\pi_0$ is faithful on $(A/M)\otimes (B/N)$, so that it is faithful on $(A/M) \widehat{\otimes} (B/N)$ \cite[Theorem 2]{jain}. Thus $q(J) =0$. Finally, since $(A\widehat{\otimes} B)/J$ is isomorphic to $(A/M) \widehat{\otimes} (B/N)$, using Theorem \ref{theo1}, it is easy to see that $M$ and $N$ are maximal in $A$ and $B$, respectively.
\end{pf}

\vspace*{2mm}
Finally, we obtain a complete picture of the lattice of closed ideals of  $\B(H)\widehat{\otimes}\\
\B(H)$.
\begin{thm}
The only non trivial closed ideals of $\B(H) \widehat\otimes \B(H)$ are $\K(H) \widehat\otimes \\ \K(H), \B(H) \widehat\otimes \K(H), \K(H) \widehat\otimes \B(H)$ and 
$\B(H) \widehat\otimes \K(H) + \K(H) \widehat\otimes \B(H)$, H being an infinite dimensional Hilbert space. 
\end{thm}

\begin{pf}
It is known that $\K(H)$ is the only non-trivial closed ideal of $\B(H)$, so using Proposition \ref{min} and Theorem \ref{max}, we have $\K(H) \oop \K(H)$ and $\B(H) \oop
\K(H) \\ + \K(H) \oop \B(H)$ as the unique minimal and maximal closed ideals of $\B(H) \oop \B(H)$, respectively. Now, consider any non trivial closed ideal $K$ of 
$\B(H) \oop \B(H)$. Using Proposition \ref{prop2}, and the fact that any proper closed ideal in a ring with unity must be contained in some maximal ideal, we notice that
\begin{equation}\label{ii}
 \K(H) \widehat \otimes \K(H) \subseteq K \subseteq \B(H) \oop \K(H) + \K(H) \oop \B(H).
\end{equation}
Let us denote $I= \B(H) \oop \K(H)$ and $ J=  \K(H) \oop \B(H)$. We first claim that 
\begin{equation}\label{1}
K \cap (I +J) = K \cap I + K\cap J.
\end{equation}
Consider any $x\in K \cap (I+ J)$. By %Corollary \ref{quasiid}, \cite[Lemma 3.3]{ass},
Lemma \ref{app}, \cite[Lemma 1.4.9]{kan}, and Proposition \ref{prop3.1}, $I +J$ has a bounded approximate identity. So using Cohen's factorization Theorem \cite{cohen}, there exist $y,z  \in (I+ J)$ such that $x = yz$, and $z$ belongs to the closed left ideal generated by $x$. Thus $z \in K$, which further gives $x \in K \cap I + K \cap J$. The other inclusion is easy.

Now $K \cap I$ and $K \cap J$ are (non-zero) closed ideals of $I$ and $J$ respectively, so using Theorem \ref{theo2}, we can write 
\begin{equation}\label{2}
 K \cap I = L \oop \K(H) \quad \text{and} \quad K \cap J= \K(H) \widehat \otimes M,
\end{equation}
where $L$ and $M$ are either $\B(H)$ or $\K(H)$. Using equations (\ref{ii}), (\ref{1}) and (\ref{2}), we have
$$K = L \oop \K(H) + \K(H) \oop M,$$
which proves the result.   

% We know that $K(H)$ is the only closed ideal of B(H), so using Proposition \ref{min} and Theorem \ref{max}, we obtain that $K(H) \widehat\otimes K(H)$ and $B(H) \widehat\otimes K(H) + K(H) \widehat\otimes B(H)$ are the minimal and maximal closed ideals of $B(H) \widehat\otimes B(H)$. Now, consider a non trivial closed ideal $K$ of $B(H) \widehat\otimes B(H)$. Using the fact that any proper closed ideal in a ring with unity must be contained in some maximal ideal, we notice that $K(H) \widehat \otimes K(H) \subset K \subset B(H) \widehat\otimes K(H) + K(H) \widehat\otimes B(H)$. Now $K \cap (B(H) \widehat\otimes K(H))$ is a closed ideal of $B(H) \widehat\otimes K(H)$, so using Theorem \ref{theo2}, we can write $$K \cap (B(H) \widehat\otimes K(H)) = L \widehat\otimes K(H),$$ where $L=B(H)$ or $K(H)$. Similarly, $K \cap (K(H) \widehat\otimes B(H))= K(H) \widehat \otimes M$, where $M=B(H)$ or $K(H)$. Also, it is easy to see that $B(H) \widehat\otimes K(H) +K(H) \widehat\otimes B(H)$ has quasi central approximate identity \cite[Lemma 3.3]{ass}, so that the closed ideal $$K \cap (B(H) \widehat\otimes K(H)) + K \cap (K(H) \widehat\otimes B(H))$$ is dense in $K \cap (B(H) \widehat\otimes K(H)+ K(H)\widehat\otimes B(H))$. Hence $K= L\widehat\otimes K(H) + K(H) \widehat\otimes M$, where $L$ and $M$ are either $K(H)$ or $B(H)$, proving the result.   

\end{pf}

%%%%%%%%%%%%%%%%%%%%%%%%%%%%%%%%%%%%%%%%%%%%%%%%%%%%%%%%%%%%%%%%%%%

\end{document}